\documentclass[12pt]{article}%
\usepackage{amsfonts}
\usepackage{amsmath}
\usepackage{amssymb}
\usepackage{graphicx}
\usepackage[small,compact]{titlesec}
\usepackage[margin=1in]{geometry}
\usepackage[title]{appendix}
\usepackage[onehalfspacing]{setspace}%
\providecommand{\U}[1]{\protect\rule{.1in}{.1in}}
\newtheorem{theorem}{Theorem}
\newtheorem{assumption}{Assumption}
\newtheorem{corollary}{Corollary}

\newtheorem{lemma}{Lemma}

\newenvironment{proof}[1][Proof]{\noindent \textbf{#1. } }{\  \rule{0.5em}{0.5em}}

\expandafter
\def \expandafter \normalsize \expandafter{\normalsize \setlength \abovedisplayskip{5pt plus 2pt minus 3pt}}
\expandafter
\def \expandafter \normalsize \expandafter{\normalsize \setlength \abovedisplayshortskip{0pt plus 2pt}}
\expandafter
\def \expandafter \normalsize \expandafter{\normalsize \setlength \belowdisplayskip{5pt plus 2pt minus 3pt}}
\expandafter
\def \expandafter \normalsize \expandafter{\normalsize \setlength \belowdisplayshortskip{2pt plus 2pt}}
\begin{document}

\title{Weak convergence to derivatives of fractional Brownian motion\thanks{We are
grateful to the editor, Peter C.\ B.\ Phillips, and two anonymous referees for
very helpful and detailed comments. We thank the Danish National Research
Foundation for financial support (DNRF Chair grant number DNRF154).}}
\author{S\o ren Johansen\thanks{Email:\ \texttt{Soren.Johansen@econ.ku.dk}}\\University of Copenhagen\\and \textsc{creates}
\and Morten \O rregaard Nielsen\thanks{Corresponding author.
Email:\ \texttt{mon@econ.au.dk}}\\Aarhus University}
\maketitle

\begin{abstract}
It is well known that, under suitable regularity conditions, the normalized
fractional process with fractional parameter $d$ converges weakly to
fractional Brownian motion for $d>1/2$. We show that, for any non-negative
integer $M$, derivatives of order $m=0,1,\dots,M$ of the normalized fractional
process with respect to the fractional parameter $d$, jointly converge weakly
to the corresponding derivatives of fractional Brownian motion. As an
illustration we apply the results to the asymptotic distribution of the score
vectors in the multifractional vector autoregressive model.




\end{abstract}

\section{Introduction}

\label{sec:intro}

The $p$-dimensional fractionally integrated process of Type~II (e.g.,
Marinucci and Robinson, 1999), is given by%
\begin{equation}
\Delta_{+}^{-d}\xi_{t}=(1-L)_{+}^{-d}\xi_{t}=\sum_{n=0}^{t-1}\pi_{n}%
(d)\xi_{t-n}=\sum_{n=1}^{t}\pi_{t-n}(d)\xi_{n},\quad t=1,2,\dots.
\label{Fractional}%
\end{equation}
This expression defines the operator $\Delta_{+}^{-d}=(1-L)_{+}^{-d}$ as a
finite sum, and the fractional coefficients $\pi_{n}(d)$ are defined by the
binomial expansion of $(1-z)^{-d}$. That is,%
\[
\pi_{n}(d)=(-1)^{n}\binom{-d}{n}=d(d+1)\dots(d+n-1)/n!\sim cn^{d-1}%
\]
with \textquotedblleft$\sim$\textquotedblright\ denoting that the ratio of the
left- and right-hand sides converges to one. The parameter $d$ is called the
memory parameter, which we assume satisfies $d>1/2$. Throughout, $\xi_{t}$ is
a $p$-dimensional linear process,%
\begin{equation}
\xi_{t}=C(L)\varepsilon_{t}=\sum_{j=-\infty}^{\infty}C_{j}\varepsilon_{t-j},
\label{Linear}%
\end{equation}
for some $p\times p$ coefficient matrices $C_{j}$ and a $p$-dimensional
innovation sequence, $\varepsilon_{t}$, which is independently and identically
distributed (i.i.d.)\ with mean zero and variance matrix $\Sigma$ (precise
conditions will be given in Section~\ref{sec:main 1}).

We define the normalized process $Z_{\left\lfloor Tr\right\rfloor
}(d)=T^{1/2-d}\Delta_{+}^{-d}\xi_{\left\lfloor Tr\right\rfloor }$ for $d>1/2$
and $r\in\lbrack0,1]$, where $\left\lfloor \cdot\right\rfloor $ denotes the
integer-part of the argument. The functional central limit
theorem\footnote{Even earlier results were available for the so-called Type~I
process; e.g.\ Davydov (1970) and Taqqu (1975).} (FCLT) for $Z_{\left\lfloor
Tr\right\rfloor }(d)$ was proved by Akonom and Gourieroux (1987) for ARMA
processes $\xi_{t}$, and by Marinucci and Robinson (2000) for linear processes
$\xi_{t}$ with coefficients satisfying a summability condition; see
Assumption~\ref{assn:xi} below. In particular, these authors showed that
\begin{equation}
Z_{\left\lfloor Tr\right\rfloor }(d)=T^{1/2-d}\sum_{n=0}^{\left\lfloor
Tr\right\rfloor -1}\pi_{n}(d)\xi_{\left\lfloor Tr\right\rfloor -n}%
\Rightarrow\Gamma(d)^{-1}\int_{0}^{r}(r-s)^{d-1}\mathsf{d}W(s)=W(r;d),
\label{AG}%
\end{equation}
where $\Gamma(\cdot)$ is the Gamma function, $W$ is Brownian motion with
variance matrix $C(1)\Sigma C(1)^{\prime}$, $C(1)=\sum_{j=-\infty}^{\infty
}C_{j}$, and \textquotedblleft$\Rightarrow$\textquotedblright\ denotes weak
convergence in the space of c\`{a}dl\`{a}g functions on $[0,1]$ endowed with
the Skorokhod topology; see Billingsley (1968) for a general treatment. That
is, the normalized process $Z_{[Tr]}(d)$ converges weakly to fractional
Brownian motion (fBm), $W(r;d)$, which is also of Type~II; see Marinucci and
Robinson (1999) for a detailed comparison of Types~I and~II fBm.

In fact, the results in Marinucci and Robinson (2000) also imply weak
convergence of the derivative of $\Delta_{+}^{-d}\xi_{t}$, suitably
normalized. We use $\mathsf{D}_{d}^{m}$ to denote the $m$'th order derivative
with respect to $d$. Differentiating term-by-term we find $\mathsf{D}_{d}%
\pi_{n}(d)=\pi_{n}(d)\sum_{k=0}^{n-1}(k+d)^{-1}$; see Appendix~A of Johansen
and Nielsen (2016) for additional details on the fractional coefficients and
their derivatives. With this notation, Marinucci and Robinson (2000) proved
that%
\begin{align}
Z_{\left\lfloor Tr\right\rfloor }^{\ast}(d)  &  =T^{1/2-d}(\log T)^{-1}%
\mathsf{D}_{d}\Delta_{+}^{-d}\xi_{\left\lfloor Tr\right\rfloor }\nonumber\\
&  =T^{1/2-d}(\log T)^{-1}\sum_{n=1}^{\left\lfloor Tr\right\rfloor -1}\pi
_{n}(d)\sum_{k=0}^{n-1}(k+d)^{-1}\xi_{\left\lfloor Tr\right\rfloor
-n}\Rightarrow W(r;d). \label{MR}%
\end{align}
Thus, because of the factor $\sum_{k=0}^{n-1}(k+d)^{-1}\sim\log n$, a
different normalization is needed, but the weak limit is still fBm.

Related to (\ref{AG}) and (\ref{MR}), Hualde (2012) showed the limit
result\footnote{There is a missing minus sign in either (6) or (8) in Hualde
(2012). Of course, this is irrelevant for the marginal distribution of
$A(r;d)$ because $A(r;d)$ is a zero-mean Gaussian process. However, the sign
is critical when considering the joint distribution of $A(r;d)$ and $W(r;d)$,
for example.}
\begin{equation}
H_{\left\lfloor Tr\right\rfloor }(d)=T^{1/2-d}\sum_{n=0}^{\left\lfloor
Tr\right\rfloor -1}\pi_{n}(d)(-\sum_{k=n}^{T}(k+d)^{-1})\xi_{\left\lfloor
Tr\right\rfloor -n}\Rightarrow A(r;d), \label{H}%
\end{equation}
where $A(r;d)=\Gamma(d)^{-1}\int_{0}^{r}\log(r-s)(r-s)^{d-1}\mathsf{d}W(s)$
was denoted a \textquotedblleft modified fBm\textquotedblright. The derivation
of (\ref{H}) was motivated by a regression analysis of so-called
\textquotedblleft unbalanced cointegration\textquotedblright, where the
process $A(r;d)$ enters in the asymptotic distribution theory; see Hualde
(2012, 2014). Note, however, that $A(r;d)=\Gamma^{-1}(d)\mathsf{D}_{d}%
(\Gamma(d)W(r;d))$ is not the derivative of fBm.

In this paper, we prove related results for weak convergence of the
derivatives with respect to $d$ of $Z_{\left\lfloor Tr\right\rfloor }(d)$ to
corresponding derivatives of fBm. Differentiating term-by-term as in
(\ref{MR}) we find%
\begin{equation}
\mathsf{D}_{d}Z_{t}(d)=\sum_{n=0}^{t-1}\mathsf{D}_{d}(T^{1/2-d}\pi_{n}%
(d))\xi_{t-n}=T^{1/2-d}\sum_{n=0}^{t-1}(-\log T+\sum_{k=0}^{n-1}\frac{1}%
{k+d})\pi_{n}(d)\xi_{t-n}. \label{m=1}%
\end{equation}
In the general case, the coefficients in the linear representation of
$\mathsf{D}_{d}^{m}Z_{t}(d)$ will be calculated by recursion; see
Section~\ref{sec:main 2}\ and Lemma~\ref{Lemma Recursion}. Note the relation%
\begin{equation}
\mathsf{D}_{d}Z_{t}(d)=(\log T)(Z_{t}^{\ast}(d)-Z_{t}(d)).
\label{samlet formel}%
\end{equation}

In recent work, Johansen and Nielsen (2021) generalize earlier work on
statistical inference in the fractionally cointegrated vector autoregressive
model (Johansen and Nielsen, 2012b) to allow each variable in the multivariate
process to have its own fractional parameter (integration order). They call
this the \textquotedblleft multifractional\textquotedblright\ vector
autoregressive model. One interpretation of this model is a generalization of
Hualde's (2014) bivariate unbalanced cointegrated regression model to a
multivariate system framework. Johansen and Nielsen (2021) show that, in this
setting, the derivative $\mathsf{D}_{d}Z_{\left\lfloor Tr\right\rfloor }(d)$
and its weak limit $\mathsf{D}_{d}W(r;d)$ play an important role in the
asymptotic distribution theory for the maximum likelihood estimators of the
fractional parameters. We present some details of this analysis in
Section~\ref{sec:app} to motivate and apply our results.

In Section~\ref{sec:main 1} we show that the result (\ref{H}) of Hualde (2012)
can be generalized to allow for weights $(-\sum_{k=n}^{T}(k+d)^{-1})^{m}$ for
any integer $m\geq0$. In Section~\ref{sec:main 2} we use this result together
with (\ref{AG}) of Marinucci and Robinson (2000) to show weak convergence of
$\mathsf{D}_{d}^{m}Z_{\left\lfloor Tr\right\rfloor }(d)$ to derivatives of
fBm. The application of our results to the multifractional cointegration model
is given in Section~\ref{sec:app}, and some concluding remarks are given in
Section~\ref{sec:conclusion}. In the next section, however, we first consider
$m=1$, because the arguments simplify substantially in that case.

\section{Weak convergence of the derivative $\mathsf{D}_{d}Z_{\left\lfloor
Tr\right\rfloor }(d)$}

\label{sec:first}

In this section, we apply the results of Marinucci and Robinson (2000) in
(\ref{AG}) and Hualde (2012) in (\ref{H}) to show that the first derivative of
the fractional process, i.e.$\ \mathsf{D}_{d}Z_{\left\lfloor Tr\right\rfloor
}(d)$, converges weakly to $\mathsf{D}_{d}W(r;d)$. Precise conditions under
which the results hold will be stated in Section~\ref{sec:main 1} before we
give the general results.

The derivative $\mathsf{D}_{d}Z_{\left\lfloor Tr\right\rfloor }(d)$ is
rewritten, using (\ref{m=1}) and $\sum_{k=0}^{n-1}(k+d)^{-1}=\sum_{k=0}%
^{T}(k+d)^{-1}-\sum_{k=n}^{T}(k+d)^{-1}$, as
\begin{align}
\mathsf{D}_{d}Z_{\left\lfloor Tr\right\rfloor }(d)  &  =\sum_{n=0}%
^{\left\lfloor Tr\right\rfloor -1}(\mathsf{D}_{d}T^{1/2-d}\pi_{n}%
(d))\xi_{\left\lfloor Tr\right\rfloor -n}\nonumber\\
&  =\sum_{n=0}^{\left\lfloor Tr\right\rfloor -1}T^{1/2-d}\pi_{n}(d)(-\log
T+\sum_{k=0}^{T}(k+d)^{-1})\xi_{\left\lfloor Tr\right\rfloor -n}\nonumber\\
&  \quad+\sum_{n=0}^{\left\lfloor Tr\right\rfloor -1}T^{1/2-d}\pi_{n}%
(d)(-\sum_{k=n}^{T}(k+d)^{-1})\xi_{\left\lfloor Tr\right\rfloor -n}\nonumber\\
&  =(-\log T+\sum_{k=0}^{T}(k+d)^{-1})Z_{\left\lfloor Tr\right\rfloor
}(d)+H_{\left\lfloor Tr\right\rfloor }(d). \label{first deriv}%
\end{align}
Here, $Z_{\left\lfloor Tr\right\rfloor }(d)\Rightarrow W(r;d)$ and
$H_{\left\lfloor Tr\right\rfloor }(d)\Rightarrow A(r;d)$ by (\ref{AG}) and
(\ref{H}), respectively. Strictly speaking, (\ref{AG}) and (\ref{H}) need to
hold jointly, but that is a consequence of Theorem~\ref{thm:Hgen} below.

To evaluate the factor $-\log T+\sum_{k=0}^{T}(k+d)^{-1}$ in
(\ref{first deriv}), recall the following definition and series expansion of
the Digamma function,%
\[
\psi(d)=\mathsf{D}_{d}\log\Gamma(d)=-\gamma-\sum_{k=0}^{\infty}((k+d)^{-1}%
-(k+1)^{-1})\text{ for }d\neq0,-1,\dots,
\]
where $\gamma=\lim_{T\rightarrow\infty}(\sum_{k=1}^{T}k^{-1}-\log
T)=0.577\dots$ is the Euler-Mascheroni constant; see Abramowitz and Stegun
(1972, eqns.\ 6.3.1 and 6.3.16). We then find that%
\begin{align}
-\log T+\sum_{k=0}^{T}\frac{1}{k+d}  &  =-(\log T-\sum_{k=1}^{T}k^{-1}%
)-(\sum_{k=0}^{T-1}(k+1)^{-1}-\sum_{k=0}^{T}(k+d)^{-1})\nonumber\\
&  \rightarrow\gamma-\sum_{k=0}^{\infty}((k+1)^{-1}-(k+d)^{-1})=-\psi(d).
\label{conv digamma}%
\end{align}

Finally we prove that%
\begin{equation}
\mathsf{D}_{d}\int_{0}^{r}(r-s)^{d-1}\mathsf{d}W(s)=\int_{0}^{r}%
\log(r-s)(r-s)^{d-1}\mathsf{d}W(s)\text{ for }d>1/2. \label{diff under}%
\end{equation}
By definition of the derivative,%
\[
\mathsf{D}_{d}\int_{0}^{r}(r-s)^{d-1}\mathsf{d}W(s)-\int_{0}^{r}%
\log(r-s)(r-s)^{d-1}\mathsf{d}W(s)=\lim_{\delta\rightarrow0}K_{d}(\delta),
\]
where%
\[
K_{d}(\delta)=\delta^{-1}\int_{0}^{r}(r-s)^{d-1}((r-s)^{\delta}-1-\delta
\log(r-s))\mathsf{d}W(s).
\]
By the mean value theorem,%
\[
(r-s)^{\delta}-1-\delta\log(r-s)=\frac{1}{2}\delta^{2}\log^{2}%
(r-s)(r-s)^{\delta^{\ast}}\text{ for }|\delta^{\ast}|\leq|\delta|.
\]
Hence we find, using the Frobenius norm $\left\Vert A\right\Vert
=(\operatorname*{tr}\{A^{\prime}A\})^{1/2}$,
\begin{align*}
\left\Vert \operatorname{Var}(K_{d}(\delta))\right\Vert  &  =\frac{1}{4}%
\delta^{2}\left\Vert \operatorname{Var}(W(1))\right\Vert \int_{0}^{r}\log
^{4}(r-s)(r-s)^{2d-2+2\delta^{\ast}}\mathsf{d}s\\
&  \leq c\delta^{2}\int_{0}^{r}\log^{4}(r-s)(r-s)^{2d-2-2|\delta|}%
\mathsf{d}s=c\delta^{2}\int_{0}^{r}(\log^{4}s)s^{|\delta|}s^{2d-2-3|\delta
|}\mathsf{d}s\\
&  \leq c\delta^{2}\int_{0}^{r}s^{2d-2-3|\delta|}\mathsf{d}s=c\delta
^{2}r^{2d-1-3|\delta|}\rightarrow0\text{ as }\delta\rightarrow0
\end{align*}
because $2d-1>0$. This proves (\ref{diff under}).

Combining these results, it follows that
\begin{align}
\mathsf{D}_{d}Z_{\left\lfloor Tr\right\rfloor }(d)  &  \Rightarrow
-\psi(d)W(r;d)+A(r;d)=\Gamma(d)^{-1}\int_{0}^{r}(-\psi(d)+\log
(r-s))(r-s)^{d-1}\mathsf{d}W(s)\nonumber\\
&  =\int_{0}^{r}\mathsf{D}_{d}(\Gamma(d)^{-1}(r-s)^{d-1})\mathsf{d}%
W(s)=\mathsf{D}_{d}W(r;d). \label{DZT}%
\end{align}
Thus, the first derivative of the fractional process $Z_{\left\lfloor
Tr\right\rfloor }(d)$ converges weakly to the first derivative of the fBm
$W(r;d)$. Interestingly, the above arguments leading to (\ref{DZT}) required
only the weak convergences in (\ref{AG}) and (\ref{H}) (jointly) together with
some well-known results regarding the Digamma function. Consequently, our
result (\ref{DZT}) holds whenever (\ref{AG}) and (\ref{H}) hold jointly. In
the next two sections we will prove the corresponding result for derivatives
of any order under precisely stated conditions.

\section{A generalization of the result of Hualde (2012)}

\label{sec:main 1}

In this section, we generalize the result (\ref{H}) of Hualde (2012). To this
end, we define the processes%
\begin{align}
H_{m,\left\lfloor Tr\right\rfloor }(d)  &  =T^{1/2-d}\sum_{n=0}^{\left\lfloor
Tr\right\rfloor -1}\pi_{n}(d)(-\sum_{k=n}^{T}(k+d)^{-1})^{m}\xi_{\left\lfloor
Tr\right\rfloor -n}\nonumber\\
&  =T^{1/2-d}\sum_{n=1}^{\left\lfloor Tr\right\rfloor }\pi_{\left\lfloor
Tr\right\rfloor -n}(d)(-\sum_{k=\left\lfloor Tr\right\rfloor -n}^{T}%
(k+d)^{-1})^{m}\xi_{n},\quad m=0,1,2,\dots, \label{Hm def}%
\end{align}
so that $Z_{\left\lfloor Tr\right\rfloor }(d)=H_{0,\left\lfloor
Tr\right\rfloor }(d)$ and $H_{\left\lfloor Tr\right\rfloor }%
(d)=H_{1,\left\lfloor Tr\right\rfloor }(d)$. In Theorem~\ref{thm:Hgen} below
we find the joint weak limit of $H_{m,\left\lfloor Tr\right\rfloor }(d)$,
$m=0,1,\dots,M$, for any non-negative integer $M$, but first we state our assumptions.

\begin{assumption}
\label{assn:xi}The $p$-dimensional process $\xi_{t}$ is such that%
\[
\xi_{t}=\sum_{j=-\infty}^{\infty}C_{j}\varepsilon_{t-j},\quad\sum
_{j=0}^{\infty}\sum_{k=j+1}^{\infty}\left(  \left\Vert C_{k}\right\Vert
^{2}+\left\Vert C_{-k}\right\Vert ^{2}\right)  <\infty,
\]
where the $C_{j}$ are $p\times p$ deterministic matrices and $C(1)=\sum
_{j=-\infty}^{\infty}C_{j}$ has full rank, $p$.
\end{assumption}

\begin{assumption}
\label{assn:eps}The $p$-dimensional process $\varepsilon_{t}$ in
Assumption~\ref{assn:xi} is i.i.d.\ with%
\[
E(\varepsilon_{t})=0,\quad E(\varepsilon_{t}\varepsilon_{t}^{\prime}%
)=\Sigma,\quad E\left\Vert \varepsilon_{t}\right\Vert ^{q}<\infty,
\]
for some $q>\max\{2,2/(2d-1)\}$, $d>1/2$, and $\Sigma$ positive definite.
\end{assumption}

We note that the moment condition in Assumption~\ref{assn:eps} is in fact
necessary; see Johansen and Nielsen (2012a). The rank conditions in
Assumptions~\ref{assn:xi}--\ref{assn:eps} ensure that the long-run variance of
$\xi_{t}$ is positive definite.

Assumptions~\ref{assn:xi}--\ref{assn:eps} are identical to the corresponding
conditions in Hualde (2012) and Marinucci and Robinson (2000). Thus,
(\ref{AG}), (\ref{H}), and the results in Section~\ref{sec:first}, and in
particular the weak convergence in (\ref{DZT}), all hold under
Assumptions~\ref{assn:xi}--\ref{assn:eps}.

\begin{theorem}
\label{thm:Hgen}Under Assumptions~\ref{assn:xi}--\ref{assn:eps} it holds that,
for $m=0,1,2,\dots$,
\begin{equation}
H_{m,\left\lfloor Tr\right\rfloor }(d)\Rightarrow A_{m}(r;d), \label{Thm1-eqn}%
\end{equation}
where $A_{m}(r;d)=\Gamma(d)^{-1}\int_{0}^{r}(\log(r-s))^{m}(r-s)^{d-1}%
\mathsf{d}W(s)$. For any non-negative integer $M$, the convergence in
(\ref{Thm1-eqn}) holds jointly for $m=0,1,\dots,M<\infty$.
\end{theorem}

\begin{proof}
The main steps of the proof are identical to those in Marinucci and Robinson
(2000) and Hualde (2012), so we focus on the relevant differences. We give the
proof for a fixed $m$. Joint convergence follows by application of the
Cram\'{e}r-Wold device and the same proof.

Marinucci and Robinson (2000) generalize the results of Einmahl (1989) to
short-range dependent variables, so they can construct copies in distribution
of $\xi_{t}$, say $\hat{\xi}_{t}$, and independent $w_{t}$ that are
i.i.d.$N(0,\Sigma)$ on the same probability space. We further define
$S_{j}=\sum_{t=1}^{j}\hat{\xi}_{t}$, $V_{j}=C(1)\sum_{t=1}^{j}w_{t}$,
$S_{0}=V_{0}=0$, and consider below the difference $S_{j}-V_{j}$, which is
possible because $S_{j}$ and $V_{j}$ are defined on the same probability
space. Specifically, based on results of Einmahl (1989, Theorems~1, 2, and~4),
Marinucci and Robinson (2000, Lemma~2) show that $\sup_{1\leq j\leq T}%
|S_{j}-V_{j}|=o_{P}(T^{1/s})$ for $2<s<q$, where $q$ is given in
Assumption~\ref{assn:eps}. As in Hualde (2012), we define
\[
\hat{H}_{m,\left\lfloor Tr\right\rfloor }(d)=T^{1/2-d}\sum_{n=1}^{\left\lfloor
Tr\right\rfloor }\pi_{\left\lfloor Tr\right\rfloor -n}(d)(-\sum
_{k=\left\lfloor Tr\right\rfloor -n}^{T}(k+d)^{-1})^{m}\hat{\xi}_{n},\quad
m=0,1,2,\dots.
\]
That is, $\hat{H}_{m,\left\lfloor Tr\right\rfloor }(d)$ is defined exactly
like $H_{m,\left\lfloor Tr\right\rfloor }(d)$ in (\ref{Hm def}), but with
$\hat{\xi}_{n}$ replacing $\xi_{n}$. Because $\hat{H}_{m,\left\lfloor
Tr\right\rfloor }(d)$ is then a copy in distribution of $H_{m,\left\lfloor
Tr\right\rfloor }(d)$, it suffices to show the required result for $\hat
{H}_{m,\left\lfloor Tr\right\rfloor }(d)$.

We then decompose $\hat{H}_{m,\left\lfloor Tr\right\rfloor }(d)=\sum_{i=1}%
^{5}Q_{iT}(r)$, where%
\begin{align*}
Q_{1T}(r)  &  =\frac{1}{\Gamma(d)}T^{-1/2}\sum_{n=1}^{\left\lfloor
Tr\right\rfloor -1}\left(  r-\frac{n}{T}\right)  ^{d-1}\left(  \log\left(
r-\frac{n}{T}\right)  \right)  ^{m}(V_{n}-V_{n-1})\mathbb{I}(\left\lfloor
Tr\right\rfloor >2),\\
Q_{2T}(r)  &  =T^{1/2-d}\sum_{n=1}^{\left\lfloor Tr\right\rfloor -1}%
\pi_{\left\lfloor Tr\right\rfloor -n}(d)(S_{n}-S_{n-1}-(V_{n}-V_{n-1}%
))(-\sum_{k=\left\lfloor Tr\right\rfloor -n}^{T}(k+d)^{-1})^{m}\mathbb{I}%
(\left\lfloor Tr\right\rfloor >2),\\
Q_{3T}(r)  &  =T^{1/2-d}\sum_{n=1}^{\left\lfloor Tr\right\rfloor -1}\left(
\pi_{\left\lfloor Tr\right\rfloor -n}(d)(-\sum_{k=\left\lfloor Tr\right\rfloor
-n}^{T}(k+d)^{-1})^{m}-\frac{(Tr-n)^{d-1}}{\Gamma(d)}\left(  \log\left(
r-\frac{n}{T}\right)  \right)  ^{m}\right) \\
&  \quad\times(V_{n}-V_{n-1})\mathbb{I}(\left\lfloor Tr\right\rfloor >2),\\
Q_{4T}(r)  &  =T^{1/2-d}(-\sum_{k=0}^{T}(k+d)^{-1})^{m}(S_{\left\lfloor
Tr\right\rfloor }-S_{\left\lfloor Tr\right\rfloor -1})\mathbb{I}(\left\lfloor
Tr\right\rfloor >2),\\
Q_{5T}(r)  &  =T^{1/2-d}\sum_{n=1}^{\left\lfloor Tr\right\rfloor }\hat{\xi
}_{n}\pi_{\left\lfloor Tr\right\rfloor -n}(d)(-\sum_{k=\left\lfloor
Tr\right\rfloor -n}^{T}(k+d)^{-1})^{m}\mathbb{I}(\left\lfloor Tr\right\rfloor
\leq2),
\end{align*}
and $\mathbb{I}(\cdot)$ denotes the indicator function. It suffices to show
that%
\begin{align}
Q_{1T}(r)  &  \Rightarrow A_{m}(r;d),\label{i1 term}\\
\sup_{0\leq r\leq1}\left\Vert Q_{iT}(r)\right\Vert  &  \overset{P}{\rightarrow
}0\text{\quad for }i=2,\ldots,5. \label{i2-5 terms}%
\end{align}
Note that the only difference between our $Q_{iT}(r)$ and the corresponding
terms in Hualde (2012), aside from notational differences, is that instead of
Hualde's $\sum_{k=\left\lfloor Tr\right\rfloor -n}^{T}(k+d)^{-1}$ and
$\log(r-n/T)$, we have $(-\sum_{k=\left\lfloor Tr\right\rfloor -n}%
^{T}(k+d)^{-1})^{m}$ and $(\log(r-n/T))^{m}$, respectively.

We first prove (\ref{i1 term}) and (\ref{i2-5 terms}) for $i=2,4,5$. These
proofs follow nearly identically to the corresponding proofs of (24) and (25)
in Hualde (2012), so we only outline the differences. First, we note that the
bound established for $m=1$ in (26) of Hualde (2012) can easily be generalized
to%
\[
\left\vert \log\left(  r-n/T\right)  \right\vert ^{m}\leq K\left(
r-n/T\right)  ^{-\alpha},\quad n=1,\ldots,\left\lfloor Tr\right\rfloor -1,
\]
for any $\alpha>0$ and some positive constant $K$ (if the bound applies for
$m=1$ and any $\alpha>0$, then clearly the bound also applies for any value of
$m$ on the left-hand side). Then the proof of (\ref{i1 term}) follows
identically to that of the corresponding term in (24) of Hualde (2012). To
prove (\ref{i2-5 terms}) for $i=2,4,5$ we can apply the same proofs as in
Hualde (2012) except with
\[
(\sum_{k=\left\lfloor Tr\right\rfloor -n}^{T}(k+d)^{-1})^{m}\leq(\sum
_{k=0}^{T}(k+d)^{-1})^{m}\leq K(\log T)^{m},
\]
where Hualde has $m=1$, and that change is inconsequential for the proofs.

It remains to prove (\ref{i2-5 terms}) for the $i=3$ term, which is the term
that involves the difference between the two factors $(-\sum_{k=n}%
^{T}(k+d)^{-1})^{m}$ and $\pi_{n}(d)$ and their corresponding limiting forms.
We bound $\sup_{0\leq r\leq1}\left\Vert Q_{3T}(r)\right\Vert $ by $\sup_{1\leq
n\leq T}\left\Vert C(1)w_{n}\right\Vert $ times%
\begin{align}
&  \sup_{0\leq r\leq1}T^{-1/2}\sum_{n=1}^{\left\lfloor Tr\right\rfloor
-1}\left\vert \frac{\pi_{n}(d)}{T^{d-1}}-\frac{1}{\Gamma(d)}\left(  \frac
{n}{T}\right)  ^{d-1}\right\vert (\sum_{k=n}^{T}(k+d)^{-1})^{m}%
\label{Q3bound1}\\
&  +\frac{1}{\Gamma(d)}\sup_{0\leq r\leq1}T^{-1/2}\sum_{n=1}^{\left\lfloor
Tr\right\rfloor -1}\left\vert (\sum_{k=n}^{T}(k+d)^{-1})^{m}-\left(  \log
\frac{n}{T}\right)  ^{m}\right\vert \left(  \frac{n}{T}\right)  ^{d-1}.
\label{Q3bound2}%
\end{align}
For $\left\lfloor Tr\right\rfloor >2$ and any $d\geq0$,
\begin{equation}
\sup_{0\leq r\leq1}\sup_{1\leq n\leq\left\lfloor Tr\right\rfloor -1}%
(\sum_{k=n}^{T}(k+d)^{-1})^{m}\leq(\sum_{k=1}^{T}(k+d)^{-1})^{m}\sim(\log
T)^{m}, \label{sumbound}%
\end{equation}
and thus the proof that $(\ref{Q3bound1})=o(1)$ is identical to that in (29)
of Hualde (2012) except the logarithmic term is raised to the power $m$, which
is inconsequential. Next, (\ref{Q3bound2}) is bounded by $\Gamma(d)^{-1}\leq
K$ times%
\begin{align}
&  \sup_{0\leq r\leq1}T^{-1/2}\sum_{n=1}^{\left\lfloor Tr\right\rfloor
-1}\left\vert (\sum_{k=n}^{T}(k+d)^{-1})^{m}-(\sum_{k=n}^{T}k^{-1}%
)^{m}\right\vert \left(  \frac{n}{T}\right)  ^{d-1}\label{Q3bound3}\\
&  +\sup_{0\leq r\leq1}T^{-1/2}\sum_{n=1}^{\left\lfloor Tr\right\rfloor
-1}\left\vert (\sum_{k=n}^{T}k^{-1})^{m}-(\int_{n}^{T}x^{-1}\mathsf{d}%
x)^{m}\right\vert \left(  \frac{n}{T}\right)  ^{d-1}. \label{Q3bound4}%
\end{align}
To bound these terms we use the identity $x^{m}-y^{m}=(x-y)\sum_{j=0}%
^{m-1}x^{j}y^{m-1-j}$ and bound the first factor as
\[
\sum_{k=n}^{T}\frac{1}{k}-\sum_{k=n}^{T}\frac{1}{k+d}=\sum_{k=n}^{T}\frac
{d}{k(k+d)}\leq d\sum_{k=n}^{T}\frac{1}{k^{2}}\leq Kn^{-1}.
\]
Using this bound together with (\ref{sumbound}), (\ref{Q3bound3}) is bounded
by%
\begin{align}
K(\log T)^{m-1}\sup_{0\leq r\leq1}T^{-3/2}\sum_{n=1}^{\left\lfloor
Tr\right\rfloor -1}\left(  \frac{n}{T}\right)  ^{d-2}  &  \leq K(\log
T)^{m-1}T^{1/2-d}\sum_{n=1}^{T}n^{d-2}\nonumber\\
&  \leq K(\log T)^{m}T^{\max\{1/2-d,-1/2\}}\rightarrow0. \label{bound 14}%
\end{align}
Similarly,
\[
\sum_{k=n}^{T}\frac{1}{k}-\int_{n}^{T}x^{-1}\mathsf{d}x\leq\sum_{k=n}%
^{T}\left(  \frac{1}{k}-\frac{1}{k+1}\right)  =\frac{1}{n}-\frac{1}{T+1}\leq
n^{-1}%
\]
and $\sup_{0\leq r\leq1}\sup_{1\leq n\leq\left\lfloor Tr\right\rfloor -1}%
(\int_{n}^{T}x^{-1}\mathsf{d}x)^{m}\sim(\log T)^{m}$, so that (\ref{Q3bound4})
is also bounded by (\ref{bound 14}).
\end{proof}

\section{Weak convergence of $\mathsf{D}_{d}^{m}Z_{\left\lfloor
Tr\right\rfloor }(d)$}

\label{sec:main 2}

We next analyze the derivatives of the fractional process $Z_{\left\lfloor
Tr\right\rfloor }(d)$ with respect to the fractional parameter $d$,
i.e.\ $\mathsf{D}_{d}^{m}Z_{\left\lfloor Tr\right\rfloor }(d)$. In terms of
the fractional coefficients and their derivatives, $\mathsf{D}_{d}%
^{m}Z_{\left\lfloor Tr\right\rfloor }(d)$ can be defined recursively as
follows. We apply logarithmic differentiation and let%
\begin{equation}
\mathsf{D}_{d}^{m}Z_{\left\lfloor Tr\right\rfloor }(d)=\sum_{n=0}%
^{\left\lfloor Tr\right\rfloor -1}\mathsf{D}_{d}^{m}(T^{1/2-d}\pi_{n}%
(d))\xi_{\left\lfloor Tr\right\rfloor -n}=\sum_{n=0}^{\left\lfloor
Tr\right\rfloor -1}T^{1/2-d}\pi_{n}(d)R_{Tn}^{(m)}(d)\xi_{\left\lfloor
Tr\right\rfloor -n}, \label{mth deriv}%
\end{equation}
where the coefficients $R_{Tn}^{(m)}(d)$ are defined by the relation
$\mathsf{D}_{d}^{m}(T^{1/2-d}\pi_{n}(d))=T^{1/2-d}\pi_{n}(d)R_{Tn}^{(m)}(d)$.
We note that
\[
\mathsf{D}_{d}^{m+1}Z_{\left\lfloor Tr\right\rfloor }(d)=\sum_{n=0}%
^{\left\lfloor Tr\right\rfloor -1}T^{1/2-d}\pi_{n}(d)(\mathsf{D}_{d}%
R_{Tn}^{(m)}(d)+R_{Tn}^{(1)}(d)R_{Tn}^{(m)}(d))\xi_{\left\lfloor
Tr\right\rfloor -n},
\]
so that the coefficients $R_{Tn}^{(m)}(d)$ must satisfy the recursion%
\begin{align}
R_{Tn}^{(1)}(d)  &  =\mathsf{D}_{d}\log(T^{1/2-d}\pi_{n}(d))=-\log
T+\sum_{k=0}^{n-1}(k+d)^{-1},\label{recursion 1a}\\
R_{Tn}^{(m+1)}(d)  &  =\mathsf{D}_{d}R_{Tn}^{(m)}(d)+R_{Tn}^{(1)}%
(d)R_{Tn}^{(m)}(d),\quad m=1,2,\dots. \label{recursion 1b}%
\end{align}
To illustrate the recursion, the next two terms of $R_{Tn}^{(m)}(d)$ are%
\begin{align*}
R_{Tn}^{(2)}(d)  &  =-\sum_{k=0}^{n-1}(k+d)^{-2}+(-\log T+\sum_{k=0}%
^{n-1}(k+d)^{-1})^{2},\\
R_{Tn}^{(3)}(d)  &  =2\sum_{k=0}^{n-1}(k+d)^{-3}-3(-\log T+\sum_{k=0}%
^{n-1}(k+d)^{-1})\sum_{k=0}^{n-1}(k+d)^{-2}+(-\log T+\sum_{k=0}^{n-1}%
(k+d)^{-1})^{3}.
\end{align*}

There is a similar recursive definition of the derivatives of fBm. We define
$R^{(m)}(d)$ by the relation $\mathsf{D}_{d}^{m}(\Gamma(d)^{-1}(r-s)^{d-1}%
)=\Gamma(d)^{-1}(r-s)^{d-1}R^{(m)}(d)$ and find
\begin{equation}
\mathsf{D}_{d}^{m}W(r;d)=\int_{0}^{r}\mathsf{D}_{d}^{m}(\Gamma(d)^{-1}%
(r-s)^{d-1})\mathsf{d}W(s)=\Gamma(d)^{-1}\int_{0}^{r}R^{(m)}(d)(r-s)^{d-1}%
\mathsf{d}W(s). \label{fBm mth deriv}%
\end{equation}
The first equality in (\ref{fBm mth deriv}) follows by the same proof as for
(\ref{diff under}). As in (\ref{recursion 1a}) and (\ref{recursion 1b}) we
find that the functions $R^{(m)}(d)$ must satisfy the recursion%
\begin{align}
R^{(1)}(d)  &  =\mathsf{D}_{d}\log(\Gamma(d)^{-1}(r-s)^{d-1})=-\psi
(d)+\log(r-s),\label{recursion 2a}\\
R^{(m+1)}(d)  &  =\mathsf{D}_{d}R^{(m)}(d)+R^{(1)}(d)R^{(m)}(d),\quad
m=1,2,\dots. \label{recursion 2b}%
\end{align}
To compare with $R_{Tn}^{(2)}(d)$ and $R_{Tn}^{(3)}(d)$, we find
\begin{align}
R^{(2)}(d)  &  =-\psi^{(1)}(d)+(-\psi(d)+\log(r-s))^{2},\nonumber\\
R^{(3)}(d)  &  =-\psi^{(2)}(d)-3(-\psi(d)+\log(r-s))\psi^{(1)}(d)+(-\psi
(d)+\log(r-s))^{3},\nonumber
\end{align}
where $\psi^{(j)}(d)=\mathsf{D}_{d}^{j}\psi(d)=\mathsf{D}_{d}^{j+1}\log
\Gamma(d)$ denotes the polygamma function; see Abramowitz and Stegun (1972,
eqn.\ 6.4.1). The recursive formulations in (\ref{recursion 1b}) and
(\ref{recursion 2b}) are clearly much more tractable than direct calculation
for larger values of $m$. We note, in particular, the strong similarity
between the terms $R_{Tn}^{(m)}(d)$ and $R^{(m)}(d)$. For example, for $m=1$
and with $n$ replaced by $\left\lfloor Tr\right\rfloor -\left\lfloor
Ts\right\rfloor $, we find that
\[
R_{T,\left\lfloor Tr\right\rfloor -\left\lfloor Ts\right\rfloor }%
^{(1)}(d)=-\log T+\sum_{k=0}^{T}(k+d)^{-1}-\sum_{k=\left\lfloor
Tr\right\rfloor -\left\lfloor Ts\right\rfloor }^{T}(k+d)^{-1}\rightarrow
-\psi(d)+\log(r-s)=R^{(1)}(d)
\]
as $T\rightarrow\infty$; c.f.\ (\ref{conv digamma}).

We next derive the solutions to the recursions.

\begin{lemma}
\label{Lemma Recursion}Let $g(d):\mathbb{R}^{+}\rightarrow\mathbb{R}$ and
assume that $\mathsf{D}^{m}g(d)$ exists for $m=1,2,\dots$ and define
$G(d)=\int_{0}^{d}g(s)\mathsf{d}s$. Define recursively the functions
$g_{m}(d)$, $m=1,2,\dots$, by $g_{0}(d)=1$ and%
\begin{equation}
g_{m+1}(d)=\mathsf{D}_{d}g_{m}(d)+g(d)g_{m}(d). \label{recursion g}%
\end{equation}
The solution $g_{m}(d)$ of (\ref{recursion g}) is given, for $m=1,2,\dots$,
by
\begin{equation}
g_{m}(d)=e^{-G(d)}\mathsf{D}_{d}^{m}e^{G(d)}=\sum_{(\ast)}c_{(\ast)}%
\prod_{i=1}^{m}\left(  \mathsf{D}_{d}^{i}G(d)\right)  ^{j_{i}}=\sum_{(\ast
)}c_{(\ast)}\prod_{i=1}^{m}\left(  \mathsf{D}_{d}^{i-1}g(d)\right)  ^{j_{i}},
\label{solution g}%
\end{equation}
where the summation $\sum_{(\ast)}$ extends over all $m$-tuples of
non-negative integers $(j_{1},\ldots,j_{m})$ that satisfy $\sum_{i=1}%
^{m}ij_{i}=m$ and where $c_{(\ast)}=m!\prod_{i=1}^{m}(j_{i}!(i!)^{j_{i}}%
)^{-1}$.
\end{lemma}

\begin{proof}
[Proof of Lemma \ref{Lemma Recursion}]The final equality in (\ref{solution g})
follows easily because $\mathsf{D}_{d}^{i}G(d)=\mathsf{D}_{d}^{i-1}g(d)$. We
multiply (\ref{recursion g}) by $e^{G(d)}$ with derivative $\mathsf{D}%
_{d}e^{G(d)}=e^{G(d)}g(d)$ and find%
\[
e^{G(d)}g_{m+1}(d)=e^{G(d)}\mathsf{D}_{d}g_{m}(d)+e^{G(d)}g(d)g_{m}%
(d)=\mathsf{D}_{d}(e^{G(d)}g_{m}(d)),\quad m=0,1,2,\dots.
\]
It follows by iteration that
\[
e^{G(d)}g_{m+1}(d)=\mathsf{D}_{d}(e^{G(d)}g_{m}(d))=\mathsf{D}_{d}%
^{2}(e^{G(d)}g_{m-1}(d))=\dots=\mathsf{D}_{d}^{m}(e^{G(d)}g(d))=\mathsf{D}%
_{d}^{m+1}(e^{G(d)}).
\]
Dividing by $e^{G(d)}$ we have proved the first equality in (\ref{solution g}%
). The next equality in (\ref{solution g}) follows from the Fa\`{a} di Bruno
formula, see Roman (1980, Theorem~2), which states that the derivatives of a
composite function $f(y)$, $y=G(d)$, are given by
\begin{align*}
\mathsf{D}_{d}^{m}f(G(d))  &  =\sum_{(\ast)}\frac{m!}{j_{1}!j_{2}!\cdots
j_{m}!}\mathsf{D}_{y}^{j_{1}+\dots+j_{m}}f(y)\prod_{i=1}^{m}\left(
\frac{\mathsf{D}_{d}^{i}G(d)}{i!}\right)  ^{j_{i}}\\
&  =\sum_{(\ast)}c_{(\ast)}\mathsf{D}_{y}^{j_{1}+\dots+j_{m}}f(y)\prod
_{i=1}^{m}\left(  \mathsf{D}_{d}^{i}G(d)\right)  ^{j_{i}}.
\end{align*}
Inserting $f(G(d))=e^{G(d)}$ and noting that $\mathsf{D}_{y}^{j_{1}%
+\dots+j_{m}}f(y)=f(y)$ we find (\ref{solution g}).
\end{proof}

\begin{corollary}
\label{Corollary Recursion}The solutions to the recursions (\ref{recursion 1a}%
)--(\ref{recursion 1b}) and (\ref{recursion 2a})--(\ref{recursion 2b}) are
given, for $m=1,2,\dots$, by%
\[
R_{T,n}^{(m)}(d)=\sum_{(\ast)}c_{(\ast)}\prod_{i=1}^{m}(\mathsf{D}_{d}%
^{i-1}R_{T,n}^{(1)}(d))^{j_{i}}\text{ and }R^{(m)}(d)=\sum_{(\ast)}c_{(\ast
)}\prod_{i=1}^{m}(\mathsf{D}_{d}^{i-1}R^{(1)}(d))^{j_{i}},
\]
respectively, where, for $i=2,3,\dots$,%
\begin{equation}
\mathsf{D}_{d}^{i-1}R_{T,n}^{(1)}(d)=(-1)^{i-1}(i-1)!\sum_{k=0}^{n-1}%
(k+d)^{-i}\text{ and }\mathsf{D}_{d}^{i-1}R^{(1)}(d)=-\psi^{(i-1)}(d).
\label{DRT}%
\end{equation}

\end{corollary}

\begin{proof}
Apply Lemma~\ref{Lemma Recursion} with initial functions $g(d)=R_{T,n}%
^{(1)}(d)=-\log T+\sum_{k=0}^{n-1}(k+d)^{-1}$ and $g(d)=R^{(1)}(d)=-\psi
(d)+\log(r-s)$, respectively. The solutions then follow from (\ref{solution g}).
\end{proof}

We are now ready to give our main result.

\begin{theorem}
\label{thm:result}Under Assumptions~\ref{assn:xi}--\ref{assn:eps} it holds
that, for $m=0,1,2,\dots$,%
\begin{equation}
\mathsf{D}_{d}^{m}Z_{\left\lfloor Tr\right\rfloor }(d)\Rightarrow
\mathsf{D}_{d}^{m}W(r;d),\nonumber
\end{equation}
where the derivatives are given in (\ref{mth deriv}) and (\ref{fBm mth deriv}%
). The convergence holds jointly for $m=0,\dots,M<\infty$.
\end{theorem}

\begin{proof}
For $m=0$ the result is given in (\ref{AG}), so we give the proof only for
$m\geq1$. Again, joint convergence follows by application of the
Cram\'{e}r-Wold device and the same proof.

We apply Corollary~\ref{Corollary Recursion} and find that, in view of
(\ref{mth deriv}) and (\ref{fBm mth deriv}), it is enough to prove (joint)
convergence for each $(j_{1},\dots,j_{m})$ where $j_{i}\geq0:$%
\begin{equation}
\sum_{n=1}^{\left\lfloor Tr\right\rfloor }\prod_{i=1}^{m}(\mathsf{D}_{d}%
^{i-1}R_{T,\left\lfloor Tr\right\rfloor -n}^{(1)}(d))^{j_{i}}T^{1/2-d}%
\pi_{\left\lfloor Tr\right\rfloor -n}(d)\xi_{n}\Rightarrow\int_{0}^{r}%
\prod_{i=1}^{m}(\mathsf{D}_{d}^{i-1}R^{(1)}(d))^{j_{i}}\mathsf{d}W.
\label{WTS}%
\end{equation}
With this result we can get the final result by taking the linear combination
$\sum_{(\ast)}c_{(\ast)}$; see Lemma~\ref{Lemma Recursion}. Thus, we start by
analyzing $(\mathsf{D}_{d}^{i-1}R_{T,\left\lfloor Tr\right\rfloor -n}%
^{(1)}(d))^{j}$ for some $j\geq1$. We consider two cases.

\emph{The case }$i=1$: We find, see (\ref{recursion 1a}) and
(\ref{conv digamma}), that
\begin{equation}
(R_{T,\left\lfloor Tr\right\rfloor -n}^{(1)})^{j}=((-\log T+\sum_{k=0}%
^{T}(k+d)^{-1})-\sum_{k=\left\lfloor Tr\right\rfloor -n}^{T}(k+d)^{-1}%
)^{j}=(-\psi(d)-\sum_{k=\left\lfloor Tr\right\rfloor -n}^{T}(k+d)^{-1}%
)^{j}+o(1). \label{R1Tn limit}%
\end{equation}

\emph{The case }$i\geq2$: Adding and subtracting appropriately, we write
$\mathsf{D}_{d}^{i-1}R_{T,\left\lfloor Tr\right\rfloor -n}^{(1)}(d)$ in
(\ref{DRT}) as%
\begin{align*}
\mathsf{D}_{d}^{i-1}R_{T,\left\lfloor Tr\right\rfloor -n}^{(1)}(d)  &
=(-1)^{i-1}(i-1)!\sum_{k=0}^{T}(k+d)^{-i}-(-1)^{i-1}(i-1)!\sum_{k=\left\lfloor
Tr\right\rfloor -n}^{T}(k+d)^{-i}\\
&  =-\psi^{(i-1)}(d)+o(1)+u_{i,n},
\end{align*}
where the convergence of the first term follows from Abramowitz and Stegun
(1972, eqn.\ 6.4.10) because $i\geq2$, and where $u_{i,n}=-(-1)^{i-1}%
(i-1)!\sum_{k=\left\lfloor Tr\right\rfloor -n}^{T}(k+d)^{-i}$ satisfies
$|u_{i,n}|\leq K(\left\lfloor Tr\right\rfloor -n)^{-i+1}\leq K(\left\lfloor
Tr\right\rfloor -n)^{-1}$ because $i\geq2$. Thus, in the analysis of
(\ref{WTS}), we can use the approximation
\begin{equation}
(\mathsf{D}_{d}^{i-1}R_{T,\left\lfloor Tr\right\rfloor -n}^{(1)}%
(d))^{j}=\left(  -\psi^{(i-1)}(d)\right)  ^{j}+o(1)+u_{i,n}\text{ for }i\geq2.
\label{DR1Tn limit}%
\end{equation}

\emph{Analysis of (\ref{WTS})}: We insert (\ref{R1Tn limit}) and
(\ref{DR1Tn limit}) into (\ref{WTS}) and find, using (\ref{AG}) and
Theorem~\ref{thm:Hgen},
\begin{align*}
&  \sum_{n=1}^{\left\lfloor Tr\right\rfloor }\prod_{i=1}^{m}(\mathsf{D}%
_{d}^{i-1}R_{T,\left\lfloor Tr\right\rfloor -n}^{(1)}(d))^{j_{i}}T^{1/2-d}%
\pi_{\left\lfloor Tr\right\rfloor -n}(d)\xi_{n}\\
&  =\sum_{n=1}^{\left\lfloor Tr\right\rfloor }((-\psi(d)-\sum_{k=\left\lfloor
Tr\right\rfloor -n}^{T}(k+d)^{-1})^{j_{1}}\prod_{i=2}^{m}(-\psi^{(i-1)}%
(d))^{j_{i}}+o(1)+u_{i,n})T^{1/2-d}\pi_{\left\lfloor Tr\right\rfloor -n}%
(d)\xi_{n}\\
&  \Rightarrow\int_{0}^{r}(-\psi(d)+\log(r-s))^{j_{1}}\prod_{i=2}^{m}%
(-\psi^{(i-1)}(d))^{j_{i}}\mathsf{d}W=\int_{0}^{r}\prod_{i=1}^{m}%
(\mathsf{D}_{d}^{i-1}R^{(1)}(d))^{j_{i}}\mathsf{d}W,
\end{align*}
see (\ref{recursion 2a}) and (\ref{DRT}). This proves (\ref{WTS}) and hence
the desired result.
\end{proof}

\section{Application to the multifractional cointegration model}

\label{sec:app}

One motivation for the results on the weak convergence of derivatives of the
fractional process comes from the analysis of the multifractional cointegrated
vector autoregressive (MFCVAR) model; see Johansen and Nielsen (2021). Let
$d=(d_{1},\dots,d_{p})^{\prime}$ be a vector of fractional parameters and let
$b$ be a scalar fractional parameter. The MFCVAR model with parameters
$\lambda=(d,b,\alpha,\beta,\Omega)$ and no lags is given by
\begin{equation}
\Lambda_{+}(d)X_{t}=-\alpha\beta^{\prime}(\Delta_{+}^{-b}-1)\Lambda
_{+}(d)X_{t}+\varepsilon_{t},\text{\quad}t=1,\dots,T,\label{model}%
\end{equation}
where the matrix differencing operator is $\Lambda_{+}(d)=\operatorname*{diag}%
(\Delta_{+}^{d_{1}},\dots,\Delta_{+}^{d_{p}})$ and $\varepsilon_{t}$ satisfies
Assumption~\ref{assn:eps}. In particular, $\varepsilon_{t}$ is i.i.d.\ with
mean zero and variance $\Omega$.

The properties of the solution to these equations can be found from the
corresponding result for the FCVAR model studied in Johansen and Nielsen
(2012b). We denote true values by subscript zero, and in particular $d_{0p}$
denotes the $p$'th element of $d_{0}$. Now, if we define $\tilde{X}_{t}$ by
$\Delta_{+}^{d_{p}}\tilde{X}_{t}=\Lambda_{+}(d)X_{t}$, then $\tilde{X}_{t}$ is
given by the equations%
\begin{equation}
\Delta_{+}^{d_{p}}\tilde{X}_{t}=-\alpha\beta^{\prime}(\Delta_{+}^{-b}%
-1)\Delta_{+}^{d_{p}}\tilde{X}_{t}+\varepsilon_{t},\quad t=1,\dots
,T.\label{FCVAR}%
\end{equation}
These equations define the FCVAR model of Johansen and Nielsen (2012b) with
scalar fractional parameters $d_{p}$ and $b$ together with $(\alpha
,\beta,\Omega)$. It follows from Theorem~2 of Johansen and Nielsen (2012b)
that the solution to (\ref{FCVAR}), for $(d_{p},b,\alpha,\beta,\Omega
)=(d_{0p},b_{0},\alpha_{0},\beta_{0},\Omega_{0})$, is%
\[
\tilde{X}_{t}=C_{0}\Delta_{+}^{-d_{p_{0}}}\varepsilon_{t}+\Delta_{+}%
^{b_{0}-d_{0p}}Y_{t},
\]
where $C_{0}=\beta_{0\bot}(\alpha_{0\bot}^{\prime}\beta_{0\bot})^{-1}%
\alpha_{0\bot}^{\prime}$ and $Y_{t}$ is a stationary linear process satisfying
Assumption~\ref{assn:xi}. Consequently, the solution to (\ref{model}), for
$\lambda=\lambda_{0}$, satisfies
\begin{equation}
\Lambda_{+}(d_{0})X_{t}=C_{0}\varepsilon_{t}+\Delta_{+}^{b_{0}}Y_{t}%
.\label{solution}%
\end{equation}
This shows that $\Delta_{+}^{d_{i0}}X_{it}$ is in general fractional of order
zero, for $i=1,\dots,p$, so that the model (\ref{model}) allows each component
of $X_{t}$ to have its own fractional order, and is therefore called
\textquotedblleft multifractional\textquotedblright. Pre-multiplying
(\ref{solution}) by $\beta_{0}^{\prime}$ shows that $\Delta_{+}^{-b_{0}}%
\beta_{0}^{\prime}\Lambda_{+}(d_{0})X_{t}$ is also fractional of order zero;
that is, some linear combinations of the processes $\{\Delta_{+}^{d_{i0}%
-b_{0}}X_{it}\}_{i=1}^{p}$ are fractional of order zero and hence $X_{t}$ is cointegrated.

We define the i.i.d.\ process $\xi_{t}=(\alpha_{0\bot}^{\prime}\beta_{0\bot
})^{-1}\alpha_{0\bot}^{\prime}\varepsilon_{t}$ such that $C_{0}\varepsilon
_{t}=\beta_{0\bot}\xi_{t}$. The three processes $Z_{t}(b_{0})$, $Z_{t}^{\ast
}(b_{0})$, and $\mathsf{D}_{b_{0}}Z_{t}(b_{0})$ are then defined in terms of
$\xi_{t}$ as in (\ref{AG}), (\ref{MR}), and (\ref{m=1}), respectively. It
follows from the above analysis that $Z_{\left\lfloor Tr\right\rfloor }%
(b_{0})$ and $Z_{\left\lfloor Tr\right\rfloor }^{\ast}(b_{0})$ converge weakly
to fractional Brownian motion $W(r;b_{0})$ and that $\mathsf{D}_{b_{0}%
}Z_{\left\lfloor Tr\right\rfloor }(b_{0})$ converges weakly to $\mathsf{D}%
_{b_{0}}W(r;b_{0})$, and that the processes converge jointly.

To simplify the subsequent analysis we assume that $\Omega=\Omega_{0}$,
$d_{p}=d_{0p}$, $\alpha=\alpha_{0}$, and $b=b_{0}>1/2$. This allows us to
focus on the parameters that give rise to \textquotedblleft
non-standard\textquotedblright\ asymptotic distributions, and in particular to
the application of $\mathsf{D}_{b_{0}}W(r;b_{0})$. Specifically, we define the
parameters $\theta=\beta_{0\bot}^{\prime}\beta$\ (or $\beta=\beta_{0}%
+\bar{\beta}_{0\bot}\theta$ with $\bar{A}=A(A^{\prime}A)^{-1}$ for any matrix
$A$ with full rank) and $\gamma_{i}=d_{i}-d_{i0}$ for $i=1,\dots,p$, such that
$\gamma_{p}=0$. With this notation we can define the residual, using
(\ref{model}) and (\ref{solution}), as%
\[
\varepsilon_{t}(\theta,\gamma)=(I_{p}-\alpha_{0}(\beta_{0}^{\prime}%
+\theta^{\prime}\bar{\beta}_{0\bot}^{\prime})(1-\Delta_{+}^{-b_{0}}%
))\Lambda_{+}(\gamma)(C_{0}\varepsilon_{t}+\Delta_{+}^{b_{0}}Y_{t}),
\]
and the Gaussian likelihood is%
\[
L_{T}(\theta,\gamma)=-\frac{1}{2}\operatorname*{tr}\{\Omega_{0}^{-1}T^{-1}%
\sum_{t=1}^{T}\varepsilon_{t}(\theta,\gamma)\varepsilon_{t}(\theta
,\gamma)^{\prime}\}=-\frac{1}{2}\operatorname*{tr}\{\Omega_{0}^{-1}%
M_{T}(\varepsilon(\theta,\gamma),\varepsilon(\theta,\gamma))\},
\]
where $M_{T}(a,b)=T^{-1}\sum_{t=1}^{T}a_{t}b_{t}^{\prime}$. We will use this
simple model to illustrate the role of the processes $Z_{t}(b_{0})$ and
$\mathsf{D}_{b_{0}}Z_{t}(b_{0})$ and their limits in the analysis of the score
functions for $\gamma$ and $\theta$ evaluated at $\lambda_{0}$.

The derivative of $\varepsilon_{t}(\theta,\gamma)$ with respect to $\theta$ at
$\lambda=\lambda_{0}$ in the direction $\partial\theta\in\mathbb{R}%
^{(p-r)\times r}$ is denoted $\mathsf{D}_{\theta}\varepsilon_{t}%
|_{\lambda=\lambda_{0}}(\partial\theta)$ and similarly for $\mathsf{D}%
_{\gamma}\varepsilon_{t}|_{\lambda=\lambda_{0}}(\partial\gamma),\partial
\gamma\in\mathbb{R}^{p}$, but with $\partial\gamma_{p}=0$ because $\gamma
_{p}=0$. We find
\begin{align*}
\mathsf{D}_{\theta}\varepsilon_{t}|_{\lambda=\lambda_{0}}(\partial\theta) &
=-\alpha_{0}(\partial\theta)^{\prime}\bar{\beta}_{0\bot}^{\prime}(1-\Delta
_{+}^{-b_{0}})(C_{0}\varepsilon_{t}+\Delta_{+}^{b_{0}}Y_{t})\\
&  \simeq\alpha_{0}(\partial\theta)^{\prime}\Delta_{+}^{-b_{0}}\xi
_{t}=T^{b_{0}-1/2}\alpha_{0}(\partial\theta)^{\prime}Z_{t}(b_{0}),\\
\mathsf{D}_{\gamma}\varepsilon_{t}|_{\lambda=\lambda_{0}}(\partial\gamma) &
=(I_{p}-\alpha_{0}\beta_{0}^{\prime}(1-\Delta_{+}^{-b_{0}}%
))\operatorname*{diag}(\partial\gamma)\mathsf{D}_{\gamma}\Lambda_{+}%
(\gamma)|_{\gamma=0}(C_{0}\varepsilon_{t}+\Delta_{+}^{b_{0}}Y_{t})\\
&  \simeq\alpha_{0}\beta_{0}^{\prime}\operatorname*{diag}(\partial
\gamma)\mathsf{D}_{\gamma}\Delta_{+}^{\gamma-b_{0}}|_{\gamma=0}\beta_{0\bot
}\xi_{t}=-\alpha_{0}\beta_{0}^{\prime}\operatorname*{diag}(\partial
\gamma)\beta_{0\bot}\mathsf{D}_{b_{0}}\Delta_{+}^{-b_{0}}\xi_{t}\\
&  =-T^{b_{0}-1/2}(\log T)\alpha_{0}\beta_{0}^{\prime}\operatorname*{diag}%
(\partial\gamma)\beta_{0\bot}Z_{t}^{\ast}(b_{0}),
\end{align*}
where $Z_{t}$ and $Z_{t}^{\ast}$ are given in (\ref{AG}) and (\ref{MR}), and
where we use `$\simeq$' to indicate that equality holds up to a stationary
process that disappears asymptotically when we normalize the nonstationary
processes. We identify the score vector $S_{T,\theta}$ for $\theta$ from
$\mathsf{D}_{\theta}L_{T}|_{\lambda=\lambda_{0}}(\partial\theta
)=(\operatorname{vec}\partial\theta)^{\prime}S_{T,\theta}$, and similarly for
$\gamma$. We then find that%
\begin{align*}
T^{-b_{0}+1}\mathsf{D}_{\theta}L_{T}|_{\lambda=\lambda_{0}}(\partial\theta) &
\simeq-\operatorname*{tr}\{\Omega_{0}^{-1}\alpha_{0}(\partial\theta)^{\prime
}T^{1/2}M_{T}(Z(b_{0}),\varepsilon)\},\\
T^{-b_{0}+1}(\log T)^{-1}\mathsf{D}_{\gamma}L_{T}|_{\lambda=\lambda_{0}%
}(\partial\gamma) &  \simeq-\operatorname*{tr}\{\Omega_{0}^{-1}\alpha_{0}%
\beta_{0}^{\prime}\operatorname*{diag}(\partial\gamma)\beta_{0\bot}%
T^{1/2}M_{T}(Z^{\ast}(b_{0}),\varepsilon)\},
\end{align*}
and, using $\operatorname*{tr}\{A^{\prime}B\}=(\operatorname{vec}A)^{\prime
}\operatorname{vec}B$, the scores are%
\begin{align*}
T^{-b_{0}+1}S_{T,\theta} &  \simeq-\operatorname{vec}(T^{1/2}M_{T}%
(Z(b_{0}),\varepsilon)\Omega_{0}^{-1}\alpha_{0}),\\
T^{-b_{0}+1}(\log T)^{-1}S_{T,\gamma} &  \simeq-B_{0}^{\prime}%
\operatorname{vec}(T^{1/2}M_{T}(Z^{\ast}(b_{0}),\varepsilon)\Omega_{0}%
^{-1}\alpha_{0}).
\end{align*}
Here we have defined the $(p-r)r\times p$ matrix $B_{0}=(\beta_{0}^{\prime
}e_{1}\otimes\beta_{0\bot}^{\prime}e_{1},\dots,\beta_{0}^{\prime}e_{p}%
\otimes\beta_{0\bot}^{\prime}e_{p})$, with $e_{i}$ denoting the $i$'th unit
vector in $\mathbb{R}^{p}$, and used the property that $\operatorname*{tr}%
\{\beta_{0}^{\prime}\operatorname*{diag}(\phi)\beta_{0\bot}M\}=\phi^{\prime
}B_{0}^{\prime}\operatorname{vec}M$; see Theorem~2 of Johansen and Nielsen
(2021). Thus, $S_{T,\theta}\in\mathbb{R}^{(p-r)r}$ and $S_{T,\gamma}%
\in\mathbb{R}^{p}$.

We note that the product moments $T^{1/2}M_{T}(Z(b_{0}),\varepsilon)$ and
$T^{1/2}M_{T}(Z^{\ast}(b_{0}),\varepsilon)$ converge jointly to their weak
limit $\int_{0}^{1}W(r;b_{0})\mathsf{d}W^{\prime}(r)$, so the scores become
linearly dependent in the limit. We therefore use the relation
(\ref{samlet formel}) to eliminate $Z_{t}^{\ast}(b_{0})=Z_{t}(b_{0})+(\log
T)^{-1}\mathsf{D}_{b_{0}}Z_{t}(b_{0})$, and the score for $\gamma$ becomes%
\[
T^{-b_{0}+1}S_{T,\gamma}\simeq-B_{0}^{\prime}\operatorname{vec}(T^{1/2}%
M_{T}((\log T)Z(b_{0})+\mathsf{D}_{b_{0}}Z(b_{0}),\varepsilon)\Omega_{0}%
^{-1}\alpha_{0}).
\]
We can now eliminate the linear dependence in the limit by defining the new
parameter%
\[
\operatorname{vec}\tilde{\theta}=\operatorname{vec}\theta+(\log T)B_{0}%
\gamma\in\mathbb{R}^{(p-r)r}%
\]
and%
\[
\tilde{\varepsilon}_{t}(\operatorname{vec}\tilde{\theta},\gamma)=\varepsilon
_{t}(\operatorname{vec}\theta,\gamma)=\varepsilon_{t}(\operatorname{vec}%
\tilde{\theta}-(\log T)B_{0}\gamma,\gamma).
\]
Then the scores and their joint limits become%
\begin{align*}
T^{-b_{0}+1}S_{T,\tilde{\theta}} &  \simeq-\operatorname{vec}(T^{1/2}%
M_{T}(Z(b_{0}),\varepsilon)\Omega_{0}^{-1}\alpha_{0})\Rightarrow
-\operatorname{vec}(\int_{0}^{1}W(r;b_{0})\mathsf{d}W^{\prime}(r)\Omega
_{0}^{-1}\alpha_{0}),\\
T^{-b_{0}+1}S_{T,\gamma} &  =(\log T)B_{0}^{\prime}\operatorname{vec}%
(T^{1/2}M_{T}(Z(b_{0}),\varepsilon)\Omega_{0}^{-1}\alpha_{0})\\
&  \quad-B_{0}^{\prime}\operatorname{vec}(T^{1/2}M_{T}((\log T)Z(b_{0}%
)+\mathsf{D}_{b_{0}}Z(b_{0}),\varepsilon)\Omega_{0}^{-1}\alpha_{0})\\
&  =-B_{0}^{\prime}\operatorname{vec}(T^{1/2}M_{T}(\mathsf{D}_{b_{0}}%
Z(b_{0}),\varepsilon)\Omega_{0}^{-1}\alpha_{0})\\
&  \Rightarrow-B_{0}^{\prime}\operatorname{vec}(\int_{0}^{1}\mathsf{D}_{b_{0}%
}W(r;b_{0})\mathsf{d}W^{\prime}(r)\Omega_{0}^{-1}\alpha_{0}).
\end{align*}
Thus, the introduction of the derivative of the fractional process and its
limit allows one to reparametrize the score to find a mixed Gaussian
asymptotic distribution, which can then be exploited to conduct inference for
some hypotheses in the MFCVAR model. For a detailed analysis we refer to
Johansen and Nielsen (2021).

\section{Concluding remarks}

\label{sec:conclusion}

Weak convergence of derivatives of fractional processes is interesting in its
own right. However, it is also likely to find application in statistical
analysis of inference problems related to multivariate fractional processes.

Hualde (2012) motivated his result (\ref{H}) with a bivariate regression
analysis of so-called \textquotedblleft unbalanced
cointegration\textquotedblright\ (see Hualde, 2014), but also anticipated that
results like (\ref{H}) may be useful in the statistical analysis of polynomial
co-fractionality (see Johansen, 2008, and Franchi, 2010).

In Section~\ref{sec:app} we presented an application of our results in
(\ref{DZT})\ and Theorem~\ref{thm:result} to the asymptotic distribution
theory for the maximum likelihood estimators of the fractional parameters in
the so-called \textquotedblleft multifractional\textquotedblright\ vector
autoregressive model of Johansen and Nielsen (2021). In this setting, the
derivative $\mathsf{D}_{d}Z_{\left\lfloor Tr\right\rfloor }(d)$ and its weak
limit $\mathsf{D}_{d}W(r;d)$ play an important role because they allow
avoiding linear dependence in the limit and because the asymptotic
distribution is expressed in terms of both $W(r,d)$ and $\mathsf{D}_{d}W(r,d)$.

\section*{References}

\begin{enumerate}
\item
\setlength{\itemsep}{1pt}
\setlength{\parskip}{0pt}
\setlength{\parsep}{0pt}%
Abramowitz, M. and Stegun, I.A. (1972). \emph{Handbook of Mathematical
Functions}, National Bureau of Standards, Washington D.C.

\item Akonom, J., and Gourieroux, C.\ (1987). A functional central limit
theorem for fractional processes. Technical report 8801, Paris: CEPREMAP.

\item Billingsley, P. (1968). \emph{Convergence of Probability Measures}%
\textit{,} New York: Wiley.

\item Davydov, Y.A. (1970). The invariance principle for stationary processes.
\emph{Theory of Probability and Its Applications} 15, 487--498.

\item Einmahl, U.\ (1989). Extensions of results of Koml\'{o}s, Major and
Tusn\'{a}dy to the multivariate case. \emph{Journal of Multivariate Analysis}
28, 20--68.

\item Franchi, M.\ (2010). A representation theory for polynomial
cofractionality in vector autoregressive models. \emph{Econometric Theory} 26, 1201--1217.

\item Hualde, J.\ (2012). Weak convergence to a modified fractional Brownian
motion. \emph{Journal of Time Series Analysis} 33, 519--529.

\item Hualde, J.\ (2014). Estimation of long-run parameters in unbalanced
cointegration. \emph{Journal of Econometrics} 178, 761--778

\item Johansen, S.\ (2008). A representation theory for a class of vector
autoregressive models for fractional processes. \emph{Econometric Theory} 24, 651--676.

\item Johansen, S.,\ and Nielsen, M.\O .\ (2012a). A necessary moment
condition for the fractional functional central limit theorem.
\emph{Econometric Theory} 28, 671--679.

\item Johansen, S.,\ and Nielsen, M.\O .\ (2012b). Likelihood inference for a
fractionally cointegrated vector autoregressive model. \emph{Econometrica} 80, 2667--2732.

\item Johansen, S., and Nielsen, M.\O .\ (2016).\ The role of initial values
in conditional sum-of-squares estimation of nonstationary fractional time
series models. \emph{Econometric Theory} 32, 1095--1139.

\item Johansen, S.,\ and Nielsen, M.\O .\ (2021). Statistical inference in the
multifractional cointegrated VAR model. In preparation, Aarhus University.

\item Marinucci, D., and Robinson, P.M.\ (1999). Alternative forms of
fractional Brownian motion. \emph{Journal of Statistical Planning and
Inference} 80, 111--122.

\item Marinucci, D., and Robinson, P.M.\ (2000). Weak convergence of
multivariate fractional processes. \emph{Stochastic Processes and their
Applications} 86, 103--120.

\item Roman, S.\ (1980). The formula of Fa\`{a} di Bruno. \emph{American
Mathematical Monthly} 87, 805--809.

\item Taqqu, M.S.\ (1975). Weak convergence to fractional Brownian motion and
to the Rosenblatt process. \emph{Zeitschrift f\"{u}r
Wahrscheinlichkeitstheorie und Verwandte Gebiete} 31, 287--302.
\end{enumerate}

\end{document}